\documentclass{article}
\usepackage[utf8]{inputenc}
\usepackage[top = 1in, bottom = 1in, left =1in, right = 1in]{geometry}
\usepackage{amsfonts,amscd,amssymb,amsmath,amsthm, mathrsfs}
\usepackage{tikz-cd}
\usepackage{xparse}
\usepackage{enumerate}
\usepackage{multirow}
\usepackage{amsthm}
\usepackage{thmtools}
\usepackage{fullpage}
\usepackage{url}
\usepackage{algorithm}
\usepackage{algpseudocode}
\usepackage{authblk}

\usepackage{marginnote}
\usepackage{verbatim}
\usepackage{mathtools}
\usepackage{enumerate}

\usetikzlibrary{arrows,matrix,positioning,fit}
\usetikzlibrary{pgfplots.groupplots}

\usepackage[all]{xy}
\xyoption{arc}

\usepackage{ctable,xparse}

\usepackage{wrapfig,subcaption}
\usepackage[colorlinks=true,citecolor=blue]{hyperref}

\graphicspath{ {images/} }

\usepackage[
    backend=bibtex,
    style=alphabetic,
    sortlocale=deDE, 
    natbib=true,
    url=false, 
    doi=true,
    eprint=false,
    isbn=false,
    giveninits=true,
    maxbibnames=10
]{biblatex}
\newcommand{\R}{\mathbb{R}}

\newcommand{\Q}{\mathbb{Q}}

\newcommand{\LP}{\operatorname{L}}

\declaretheorem{theorem}
\declaretheorem[sibling=theorem]{lemma}

\declaretheorem[sibling=theorem]{corollary}

\newcommand{\prob}{\mathbb{P}}

\renewcommand{\Pr}{\prob}
\DeclareDocumentCommand \one { o }
{%
\IfNoValueTF {#1}
{\mathbf{1}  }
{\mathbf{1}\{{#1}\} }%
}

\newcommand{\rank}{\operatorname{rank}}

\DeclareDocumentCommand{\Prto} {o} {
\IfNoValueTF {#1}
 {\overset{\Pr}{\longrightarrow}}
 { \xrightarrow[ #1 \to \infty]{\Pr }}
}
\DeclareDocumentCommand{\Asto} {o} {
\IfNoValueTF {#1}
 {\overset{\operatorname{a.s.}}{\longrightarrow}}
 {
 \xrightarrow[ #1 \to \infty]{\operatorname{a.s.} }
 }
}
\DeclareDocumentCommand{\Mgfto} {o} {
\IfNoValueTF {#1}
{\overset{\operatorname{mgf}}{\longrightarrow}}
{ \xrightarrow[ #1 \to \infty]{\operatorname{mgf} }}
}

\DeclareDocumentCommand{\Wkto} {o} {
\IfNoValueTF {#1}
 {\overset{(d)}{\longrightarrow}}
 { \xrightarrow[ #1 \to \infty]{(d) }}
}

\DeclareDocumentCommand \LPto { O{1} }
{\overset{\operatorname{\LP^{#1}}}{\longrightarrow}}

\title{A note on large torsion in $\Q$--acyclic complexes} 
\author{Andrew Vander Werf}

\affil{Division of Applied Mathematics, Brown University, Providence, RI}

\date{\today}

\begin{document}
\maketitle

\begin{abstract}
New upper bounds on the size of the torsion group of a $\Q$--acyclic simplicial complex are introduced which depend only on the vertex degree sequence of the complex and its dimension. 
\end{abstract}

\section{Introduction}

If a tree on $n$ vertices is an edge--maximally acyclic graph on $[n]:=\{1,2,\dots,n\}$, a $k$--dimensional simplicial analogue of a tree should, at the very least, be a $k$--dimensional simplicial complex on $[n]$ with as many faces as possible while maintaining a trivial $k$th homology group. Since the top ($k$th) homology group cannot gain elements by the addition of faces of dimension less than $k$, this face--maximality condition implies a full $(k-1)$--skeleton, which has the happy side--effect of rendering all of the lower--dimensional (rational, reduced) homology groups trivial. Just as a tree on $[n]$ always has exactly $n-1$ edges by virtue of being maximally acyclic, the maximality condition we've imposed causes this $k$--dimensional analogue of a tree to have, as it turns out, exactly ${{n-1}\choose k}$ $k$--faces. 

We mention all of this to help motivate the following definition due to Kalai \cite{Kalai}. A $k$--dimensional simplicial complex $T$ on $[n]$ is called \emph{$\Q$--acyclic} if 
\begin{itemize}
    \item $T$ has a full $(k-1)$--skeleton, 
    \item $T$ has ${{n-1}\choose k}$ $k$--faces, 
    \item and $H_k(T)=0$. 
\end{itemize}
We denote the set of such complexes by $\mathscr{T}_{n,k}$, as they are meant to be interpreted as higher--dimensional analogues of combinatorial trees. Due to the presence of a full $(k-1)$--skeleton, a $\Q$--acyclic complex is defined entirely by its set of $k$--dimensional faces. We therefore identify a $\Q$--acyclic complex with its set of $k$--dimensional faces. In doing so, we denote by $|T|$ the number of top--dimensional faces in a $\Q$--acyclic complex. 

This choice of definition has been vindicated by way of numerous elegant generalizations of famous results from the classical theory of trees \cite{Kalai}, \cite{Adin} \cite{DKM1}, \cite{DKM2}, \cite{BK}, \cite{DKM3}. Additionally, this definition introduces a new aspect to the theory of trees that is only apparent in higher dimensions: torsion. We mentioned that a $\Q$--acyclic complex always has trivial \emph{rational} reduced homology in every dimension. Nearly the same can be said when we switch to integer coefficients, the only difference being that $\tilde{H}_{d-1}(T)$ (called the \emph{torsion group} of $T$) may be a nontrivial but notably finite abelian group. 

One of the most interesting properties of $\Q$--acyclic complexes is their tendency to have exceptionally large torsion groups. In particular, for given integers $n>k\geq0$, the maximum possible size of the torsion (sub)group of a $k$--dimensional simplicial complex on $[n]$ can always be achieved by a $\Q$--acyclic complex. The purpose of this note is to introduce some new upper bounds on the maximum possible size of the torsion group of a $k$--dimensional $\Q$--acyclic simplicial complex $T$ given its vertex degree sequence $\big(d_i(T):i\in[n]\big)$ where $d_i(T)$ denotes the number of top--dimensional faces of a simplicial complex which contain the vertex $i$. We will primarily build off of the following generating function result. 
\begin{theorem}[\cite{Kalai}, Theorem 3']\label{genfunc}
Let $\{e_i:i\in[n]\}$ be commuting formal variables. Then 
$$\sum_{T\in\mathscr{T}_{n,k}}|\tilde{H}_{k-1}(T)|^2\prod_{i\in[n]}e_i^{d_i(T)}=\bigg(\sum_{i\in[n]}e_i\bigg)^{{n-2}\choose k}\prod_{i\in[n]}e_i^{{n-2}\choose{k-1}}.$$
Note the special case $e_1=e_2=\cdots=e_n=k=1$ which is Caley's formula. 
\end{theorem}

For readability, we will from now on set $m_1:={{n-2}\choose{k-1}}$, $m_2:={{n-2}\choose k}$, and $m_3:={{n-1}\choose k}$, noting that 
$$m_1+m_2=m_3\quad\text{and}\quad nm_1+m_2=(k+1)m_3.$$ 
Two things we can easily observe from this generating function which will be helpful to keep in mind are that, for all $T\in\mathscr{T}_{n,k}$, we have $\sum_{i\in[n]}d_i(T)=nm_1+m_2=(k+1)m_3$ and $m_1\leq d_i(T)\leq m_3$ for all $i\in[n]$. Having established these things, the following results are presented. 

\begin{theorem}\label{Main}
For $T\in\mathscr{T}_{n,k}$ with degree sequence $(d_i:i\in[n])$, we have 
$$|\tilde{H}_{k-1}(T)|^2\leq\inf_{x\in\R^n\setminus\{\mathbf{0}\}}\frac{\prod_{i\in[n]}(x_i^2)^{{m_1}-d_i}}{\left(\sum_{i\in[n]}x_i^2\right)^{m_1}}\prod_{\tau\in T}\sum_{i\in\tau}x_i^2.$$ 
\end{theorem}

\begin{corollary}\label{oned}
For $T\in\mathscr{T}_{n,k}$, let $d$ be the degree of any vertex in $T$. Then 
$$|\tilde{H}_{k-1}(T)|^2\leq
\begin{cases}
\left(\frac{d-m_1}{m_2}\right)^{m_1}\left(1+\frac{k}{n-1}\frac{m_3-d}{d-m_1}\right)^d(k+1)^{m_3-d},& d>m_1 \\ 
\left(\frac{k}{n-1}\right)^{m_1}(k+1)^{m_2},& d=m_1
\end{cases}.$$ 
\end{corollary}

\begin{corollary}\label{simplebound}
With $T$ as in Theorem \ref{Main} and $0^0$ understood to be $1$, we have 
$$|\tilde{H}_{k-1}(T)|^2\leq\frac{\prod_{i\in[n]}(d_i-m_1)^{m_1-d_i}}{m_2^{m_1}}\left(\frac{\sum_{i\in[n]}d_i^2}{m_3}-(k+1)m_1\right)^{m_3}.$$
\end{corollary}

\begin{corollary}\label{tighterbound}
With $T$ as in Theorem \ref{Main}, suppose by relabeling vertices that $\big(d_i:i\in[n]\big)$ is in increasing order, let $i_*:=\min\{i\in[n]:d_i>m_1\}$, and let $\alpha\in\left[m_3^{-1}d_{i_*},m_3^{-1}d_n\right]\subset\left(\frac{k}{n-1},1\right]$ be the unique solution to 
$\sum_{i\in[n]}\frac{d_i-m_1}{d_i-m_1\alpha}=\alpha^{-1}$. 
Then 
$$|\tilde{H}_{k-1}(T)|^2\leq\alpha^{m_1}\prod_{i=i_*}^n\left(1+\frac{(1-\alpha)m_1}{d_i-m_1}\right)^{d_i-m_1}\leq\left(\alpha e^{(1-\alpha)(n-i_*+1)}\right)^{m_1}\leq e^{(1-\alpha)(n-i_*)m_1}.$$    
\end{corollary}

As we can see particularly clearly through Corollary \ref{oned}, the presence of a vertex with very large degree (that is, close to $m_3$) impedes large torsion. Since the sum of the degree sequence is a fixed value, this suggests that the most spectacularly sized torsion groups may be observed among those $\Q$--acyclic complexes with highly uniform degree sequences. 

\section{Linear algebra and determinantal measure}

We will in this section give only a minimal description of the objects necessary for our proofs. Please see Sections 2.1 and 3 of \cite{VanderWerf2024} for a more detailed discussion with notation that is consistent with what follows. We consider the $k$--dimensional boundary matrix $\partial$ with rows and columns indexed respectively by ${[n]\choose k}$ and ${[n]\choose{k+1}}$---where ${[n]\choose j}$ is being used here to denote the set of ordered subsets of $[n]$ of size $j$---with entries defined for each $\sigma=\{\sigma_0<\sigma_1<\cdots<\sigma_{k-1}\}\in{[n]\choose k}$ and $\tau=\{\tau_0<\tau_1<\cdots<\tau_k\}\in{[n]\choose{k+1}}$ by  
$$\partial(\sigma,\tau)=
\begin{cases}
(-1)^m, &\sigma=\tau\setminus\{\tau_m\} \\
0,&\text{otherwise}
\end{cases}.$$

From this matrix we define the submatrix $\widehat\partial$ of $\partial$ by deleting the rows of $\partial$ which contain the vertex $n$. Given a matrix $M$ with entries indexed over the set $S\times T$, for $A\subseteq S$ and $B\subseteq T$, we write $M_{A,B}$ to denote the submatrix of $M$ with rows indexed by $A$ and columns indexed by $B$. We also write $\bullet$ in either the row or column subscript to indicate that the full index set is being used. As a point of clarification for this notation, transposes are handled by the convention of writing $M_{A,B}^t$ to mean $(M_{A,B})^t=(M^t)_{B,A}$. 

It is known \cite{Kalai} that, for $T\in\mathscr{T}_{n,k}$, the submatrix $\widehat\partial_{\bullet,T}$ is square, and the modulus of its determinant is equal to $|\tilde H_{k-1}(T)|$. Kalai used this fact and the Cauchy--Binet formula to show that, for any ${[n]\choose{k+1}}\times{[n]\choose{k+1}}$ diagonal matrix $X$ with nonzero entries from a field,  
    $$\det(\widehat\partial X^2{\widehat\partial}^t)=\sum_{T\in{{[n]\choose {k+1}}\choose{{n-1}\choose k}}}\det(\widehat\partial X)_{\bullet,T}^2=\sum_{T\in\mathscr{T}_{n,k}}|\tilde{H}_{k-1}(T)|^2\det X_{T,T}^2.$$
This last fact gives us a family of probability measures $\nu=\nu^{X}$ on $\mathscr{T}_{n,k}$---each of which is identified with a probability measure on ${[n]\choose {k+1}}\choose{{n-1}\choose k}$---parameterized by $X$. Writing out the entries of this matrix explicitly, for $\sigma,\sigma'\in{[n-1]\choose k}$ we have 
$$\widehat\partial X^2\widehat\partial^t(\sigma,\sigma')=
\begin{cases}
\sum_{\tau\supset\sigma}X_\tau^2&\text{if }\sigma=\sigma'\\
\widehat\partial(\sigma,\sigma\cup\sigma')\widehat\partial(\sigma',\sigma\cup\sigma')X_{\sigma\cup\sigma'}^2,&\text{if }|\sigma\cup\sigma'|=k+1 \\
0,&\text{otherwise}
\end{cases}.$$

We will consider the following special case, also due to Kalai. For any simplex $\tau$ of dimension $k$, set $X_k(\tau,\tau):=\prod_{i\in\tau}x_i$ where $\{x_i\}_{i\in[n]}$ are nonzero elements of a field. Then 
$$(\widehat\partial X_k^2\widehat\partial^t)(\sigma,\sigma')=\begin{cases}
(\sum_{i\notin\sigma}x_i^2)\prod_{i\in\sigma}x_i^2,&\text{if }\sigma=\sigma'\\
\widehat\partial(\sigma,\sigma\cup\sigma')\widehat\partial(\sigma',\sigma\cup\sigma')\prod_{i\in\sigma\cup\sigma'}x_i^2,&\text{if }|\sigma\cup\sigma'|=k+1 \\
0,&\text{otherwise}
\end{cases}.
$$
In Kalai's proof of Theorem \ref{genfunc}, it is shown that  
\begin{equation}\label{normalizationfactor}
\det(\widehat\partial X_k^2\widehat\partial^t)=z^{m_2}\prod_{i=1}^n(x_i^2)^{m_1} 
\end{equation} 
where $z:=\sum_{i\in[n]}x_i^2$.
Thus we have a family of probability measures on $\mathscr{T}_{n,k}$ which we will denote by $\nu^x_{n,k}$. Specifically, 
$$\nu_{n,k}^x(T):=\frac{|\tilde{H}_{k-1}(T)|^2\det(X_k)_{T,T}^2}{z^{m_2}\prod_{i=1}^n(x_i^2)^{m_1}}=\frac{|\tilde{H}_{k-1}(T)|^2\prod_{i=1}^n(x_i^2)^{d_i(T)-m_1}}{z^{m_2}}.$$


Probability measures like $\nu_{n,k}^x$ which take the form $T\mapsto\frac{\det M_{\bullet,T}^2}{\det(MM^t)}$ for some matrix $M$ are called \emph{determinantal} \cite{Lyons}, \cite{Lyons2}. The following is a well known fact about determinantal probability measures (see for example Theorem 5.1 of \cite{Lyons}). 

\begin{lemma}
Let $\mu$ be the determinantal probability measure on the set $S$ defined by 
$$\mu(T):=\frac{\det M_{\bullet,T}^2}{\det(MM^t)}$$
for all $T\subseteq S$ of size $\rank M$. Let $P:=M^t(MM^t)^{-1}M$. Then, for any $B\subseteq S$ we have 
$$\mu(\{T\subseteq S:T\supseteq B\})=\det P_{B,B}
.$$
\end{lemma}

As we see from the previous lemma, for any $B\subseteq{[n]\choose{k+1}}$, we have 
$$\nu^x_{n,k}\left(\left\{T\in\mathscr{T}_{n,k}:T\supseteq B\right\}\right)=\det P^x_{B,B}
$$
where $P^x:=(\widehat\partial X_k)^t(\widehat\partial  X_k^2\widehat\partial^t)^{-1}\widehat\partial X_k$
is the orthogonal projection onto the rowspace of $\widehat\partial X_k$. Since this matrix is difficult to calculate, we note that this $P^x$ is uniquely characterized by having all of its columns in the rowspace of $\widehat\partial X_k$ and satisfying $\widehat\partial X_k(I-P^x)=0$.
M{\'e}sz{\'a}ros (\cite{Meszaros}, Lemma 14) used this characterization to determine a simple expression for $P^x$ in the case that $x$ is the all--ones vector. Pleasantly, it turned out that $$P^{(1,1,\dots,1)}=n^{-1}\partial^t\partial.$$

\begin{lemma}\label{kernel}
Let $X_k$ and $z$ be as in (\ref{normalizationfactor}). Then 
$$P^x=z^{-1}(X_k\partial^tX_{k-1}^{-2}\partial X_k).$$
In particular, for 
$\tau,\tau'\in{[n]\choose{k+1}}$ we have 
$$P^x(\tau,\tau)=z^{-1}\sum_{i\in\tau}x_i^2.$$ 
\end{lemma}
\begin{proof}
The proof is largely the same as the proof for Lemma 14 of \cite{Meszaros}, but we'll spell out the necessary additional details here for completeness. For readability, let $X=X_k$ and $Y=X_{k-1}^{-1}$. We are trying to show for all $\tau,\tau'\in{[n]\choose{k+1}}$ that 
$$P^x(\tau,\tau')=\begin{cases}z^{-1}X_\tau^2\sum_{\sigma\subset\tau}Y^2_\sigma, &\text{if }\tau=\tau' \\ 
z^{-1}\partial(\tau\cap\tau',\tau)\partial(\tau\cap\tau',\tau')X_\tau X_{\tau'}Y_{\tau\cap\tau'}^2,& \text{if }|\tau\cap\tau'|=k \\ 
0, &\text{otherwise}
\end{cases}.$$

As mentioned a moment ago, there are two things to show: that each of the columns of $P^x$ are in the rowspace of $\widehat\partial X$, and that $\widehat\partial X(I-P^x)=0$. The first statement follows immediately from the definition of $P^x$ and the standard fact that the rows of $\widehat\partial$ span the rowspace of $\partial$. As for the second statement, we have for any $\sigma\in{[n-1]\choose k}$ and $\tau\in{[n]\choose{k+1}}$ that
$$\big(\widehat\partial X(I-P^x)\big)(\sigma,\tau)=X_\tau \left(\partial(\sigma,\tau)\bigg(1-z^{-1} X_\tau^2\sum_{\sigma'\subset\tau}Y_{\sigma'}^2\bigg)-z^{-1}\sum_{\tau\neq\tau'\in{[n]\choose {k+1}}}\partial(\sigma,\tau')\partial(\tau\cap\tau',\tau)\partial(\tau\cap\tau',\tau')X_{\tau'}^2Y_{\tau\cap\tau'}^2\right).$$
So, naturally, we can ignore the $X_\tau$ on the outside and simply show that 
$$\partial(\sigma,\tau)\left(1-z^{-1} X_\tau^2\sum_{\sigma'\subset\tau}Y_{\sigma'}^2\right)- z^{-1}\sum_{\tau\neq\tau'\in{[n]\choose {k+1}}}\partial(\sigma,\tau')\partial(\tau\cap\tau',\tau)\partial(\tau\cap\tau',\tau')X_{\tau'}^2Y_{\tau\cap\tau'}^2=0$$ for all $\sigma\in{[n-1]\choose k}$ and $\tau\in{[n]\choose{k+1}}$. There are three cases to consider which depend on the size of $|\sigma\cap\tau|$. The first case is simple. If $|\sigma\cap\tau|<k-1$, then all terms are zero due to the support of $\partial$. 

For the second case, set $S:=\sigma\cap\tau$. If $|S|=k-1$, we have $r_1,r_2\in \tau\setminus\sigma$ and $r_3\in\sigma\setminus\tau$ so that $\sigma=S\cup\{r_3\}$, $\tau=S\cup\{r_1,r_2\}$, and $\partial(\sigma,\tau')\partial(\tau\cap\tau',\tau)\partial(\tau\cap\tau',\tau')$ is nonzero only if either $\tau'=S\cup\{r_1,r_3\}$ or $\tau'=S\cup\{r_2,r_3\}$. This last fact holds because $\tau'$ has to contain $\sigma=S\cup\{r_3\}$ for the first factor to be nonzero, but it also must contain one of $r_1,r_2\in\tau\setminus\sigma$ so that $|\tau\cap\tau'|=k$ which allows the second two factors to be nonzero. So we have 
\begin{align*}
&\partial(\sigma,\tau)\left(1- z^{-1} X_\tau^2\sum_{\sigma'\subset\tau}Y_{\sigma'}^2\right)- z^{-1}\sum_{\tau\neq\tau'\in{[n]\choose {k+1}}}\partial(\sigma,\tau')\partial(\tau\cap\tau',\tau)\partial(\tau\cap\tau',\tau')X_{\tau'}^2Y_{\tau\cap\tau'}^2 \\
=&-z^{-1}\sum_{\tau\neq\tau'\in{[n]\choose {k+1}}}\partial(S\cup\{r_3\},\tau')\partial((S\cup\{r_1,r_2\})\cap\tau',S\cup\{r_1,r_2\})\partial((S\cup\{r_1,r_2\})\cap\tau',\tau')X_{\tau'}^2Y_{(S\cup\{r_1,r_2\})\cap\tau'}^2 \\
=&-z^{-1}\sum_{i=1}^2\partial(S\cup\{r_3\},S\cup\{r_i,r_3\})\partial(S\cup\{r_i\},S\cup\{r_i,r_2\})\partial(S\cup\{r_i\},S\cup\{r_i,r_3\})X_{S\cup\{r_i,r_3\}}^2Y_{S\cup\{r_i\}}^2 \\
=&-z^{-1} x_{r_3}^2\sum_{i=1}^2\partial(S\cup\{r_3\},S\cup\{r_i,r_3\})\partial(S\cup\{r_i\},S\cup\{r_i,r_2\})\partial(S\cup\{r_i\},S\cup\{r_i,r_3\})=0
\end{align*}
where we can conclude that this last line is 0 because we know, thanks to M{\'e}sz{\'a}ros (\cite{Meszaros}, Lemma 13), that it would be 0 in the special case $x_{r_3}=1$.

For the third and final case to consider, if $\sigma\subset\tau$ we have 
\begin{align*}
&\partial(\sigma,\tau)\left(1-z^{-1} X_\tau^2\sum_{\sigma'\subset\tau}Y_{\sigma'}^2\right)- z^{-1}\sum_{\tau\neq\tau'\in{[n]\choose {k+1}}}\partial(\sigma,\tau')\partial(\tau\cap\tau',\tau)\partial(\tau\cap\tau',\tau')X_{\tau'}^2Y_{\tau\cap\tau'}^2 \\
=&\partial(\sigma,\tau)\left(1-z^{-1} X_\tau^2\sum_{\sigma'\subset\tau}Y_{\sigma'}^2\right)- z^{-1}\sum_{a\in[n]\setminus\tau}\partial(\sigma,\sigma\cup\{a\})\partial(\tau\cap(\sigma\cup\{a\}),\tau)\partial(\tau\cap(\sigma\cup\{a\}),\sigma\cup\{a\})X_{\sigma\cup\{a\}}^2Y_{\tau\cap(\sigma\cup\{a\})}^2 \\
=&\partial(\sigma,\tau)\left(1-z^{-1} X_\tau^2\sum_{\sigma'\subset\tau}Y_{\sigma'}^2-z^{-1} Y_\sigma^2\sum_{a\in[n]\setminus\tau}\partial(\sigma,\sigma\cup\{a\})X_{\sigma\cup\{a\}}^2\right) \\
=& z^{-1}\partial(\sigma,\tau)\left(z-X_\tau^2\sum_{\sigma'\subset\tau}Y_{\sigma'}^2-Y_\sigma^2\sum_{\tau\neq\tau'\supset\sigma}X_{\tau'}^2\right) \\ 
=&z^{-1}\partial(\sigma,\tau)\left(\sum_{i\in[n]}x_i^2-\sum_{i\in\tau}x_i^2-\sum_{i\notin\tau}x_i^2\right)=0.
\end{align*}
The proof is now complete. 
\end{proof}

\section{Proofs of the inequalities}
\begin{proof}[Proof of Theorem \ref{Main}]    
By Hadamard's inequality and the previous two lemmas, 
$$\frac{|\tilde{H}_{k-1}(T)|^2\prod_{i=1}^n(x_i^2)^{d_i(T)-m_1}}{z^{m_2}}=\nu_{n,k}^x(T)=\det P^x_{T,T}\leq\prod_{\tau\in T}\frac{\sum_{i\in\tau}x_i^2}{z}.$$
This in particular gives us an upper bound on the size of the torsion group of a given $T$, 
\begin{equation}\label{thebound}
|\tilde{H}_{k-1}(T)|^2\leq\frac{\prod_{i\in[n]}(x_i^2)^{{m_1}-d_i(T)}}{z^{m_1}}\prod_{\tau\in T}\sum_{i\in\tau}x_i^2. 
\end{equation}
This inequality holds for all $x\in(\R\setminus\{0\})^n$. By taking some but not all of the coordinates arbitrarily close to 0, the result follows. 
\end{proof}
By plugging in 1 for each $x_i$ above, we get 
$$|\tilde{H}_{k-1}(T)|^2\leq n^{-{m_1}}(k+1)^{m_3}=\left(\frac{k+1}{n}\right)^{m_1}(k+1)^{m_2},$$
which is only a slight improvement to Kalai's initial bound of $(k+1)^{m_2}$. 
We can get a more nuanced bound in terms of one of the degrees of $T$ by setting all but one $x_i$ equal to 1. 
\begin{proof}[Proof of Corollary \ref{oned}] 
Let $x_1=\sqrt{\frac{n-1}{t}}$ for some $t>0$ and set the rest to be 1. Letting $d=d_1(T)$ without loss of generality, expression (\ref{thebound}) yields 
\begin{align*}
|\tilde{H}_{k-1}(T)|^2&\leq\frac{\left(\frac{n-1}{t}\right)^{{m_1}-d}\left(k+\frac{n-1}{t}\right)^d}{\big((n-1)(1+t^{-1})\big)^{m_1}}(k+1)^{{m_3}-d}=\frac{\left(1+\frac{k}{n-1}t\right)^d}{(1+t)^{m_1}}(k+1)^{{m_3}-d}.
\end{align*}
We now want to minimize $f(t):=\left(1+\frac{k}{n-1}t\right)^d(1+t)^{-m_1}$. Its derivative is 
$$f'(t)=\left(1+\frac{k}{n-1}t\right)^{d-1}(1+t)^{-{m_1}-1}\frac{k}{n-1}\left(\left(d-{m_1}\right)t-\left({m_3}-d\right)\right).$$
If $d=m_1$, $f(t)$ is a non--increasing function, and taking $t\to\infty$ gives the desired bound in this case. Otherwise, we reach a global minimum at $t=\frac{m_3-d}{d-m_1}=\frac{m_2}{d-m_1}-1$. Plugging this into $f$ gives us 
\begin{align*}
f\left(\frac{m_3-d}{d-m_1}\right)&=\left(1+\frac{k}{n-1}\frac{m_3-d}{d-m_1}\right)^d\left(\frac{m_2}{d-m_1}\right)^{-m_1} 
\end{align*}
as desired. 
\end{proof}

Attempting to optimize more than one variable in $(\ref{thebound})$ requires knowledge of the degrees of higher--dimensional faces. We can get around this by applying the AMGM inequality to (\ref{thebound}) to get \begin{equation}\label{lesserbound}
|\tilde{H}_{k-1}(T)|^2\leq\frac{\prod_{i\in[n]}(x_i^2)^{{m_1}-d_i(T)}}{z^{m_1}}\left(\frac{\sum_{i\in[n]}d_i(T)x_i^2}{m_3}\right)^{m_3}. 
\end{equation}
One immediate but sub--optimal bound that we can gather from this is the bound found in Corollary \ref{simplebound}, 
\begin{proof}[Proof of Corollary \ref{simplebound}]
Take $x_i^2=d_i(T)-m_1$ for each $i$ in (\ref{lesserbound}).
\end{proof}
Note that in the cases where $d_i(T)=m_1$, we are really considering the limit as $x_i$ approaches $0$. Thus $\big(d_i(T)-m_1\big)^{m_1-d_i(T)}=1$ in such cases. 
We now find the true global minimum of the right hand side of (\ref{lesserbound}). 
\begin{proof}[Proof of Corollary \ref{tighterbound}]
Let 
$$f(x):=\frac{\prod_{i\in[n]}(x_i^2)^{{m_1}-d_i}}{\left(\sum_{i\in[n]}x_i^2\right)^{m_1}}\left(\frac{\sum_{i\in[n]}d_ix_i^2}{m_3}\right)^{m_3}.$$ 
By relabeling vertices, we can assume $d_1\leq d_2\leq\cdots\leq d_n$. The first thing we can notice about $f$ is that, for every $i$ with $d_i=m_1$, the best course of action is to allow $x_i$ to tend to 0 since $m_3>m_1$, but otherwise allowing $x_i\to0$ will result in explosion. 
Thus we can assume $x_i=0$ if and only if $d_i=m_1$. 

We will first consider the problem under the constraint that 
$$\log\frac{\sum_{i\in[n]}d_ix_i^2}{\sum_{i\in[n]}x_i^2}=\log\alpha$$
for some $\alpha\in[d_{i_*},d_n]$ where $i_*:=\min\{i:\in[n]:d_i>m_1\}$. If $x$ is a constrained critical point of $\log f$, there will be a Lagrange multiplier $\lambda$ such that 
$$-2(d_i-m_1)x_i^{-1}-\frac{2m_1x_i}{\sum_{j\in[n]}x_j^2}+\frac{2m_3d_ix_i}{\sum_{j\in[n]}d_jx_j^2}=\frac{2\lambda d_ix_i}{\sum_{j\in[n]}d_jx_j^2}-\frac{2\lambda x_i}{\sum_{j\in[n]}x_j^2}$$
for each $i\geq i_*$. Equivalently, 
$$\frac{x_i^2}{\sum_{j\in[n]}d_jx_j^2}=\frac{d_i-m_1}{m_3d_i-m_1\alpha-(d_i-\alpha)\lambda}$$
for each $i\geq i_*$. In order for $x$ to be in the constraint set, $\lambda$ therefore needs to satisfy both  
$$\sum_{i\in[n]}\frac{d_i-m_1}{m_3d_i-m_1\alpha-(d_i-\alpha)\lambda}=\alpha^{-1}\quad\text{and}\quad\sum_{i\in[n]}\frac{(d_i-m_1)d_i}{m_3d_i-m_1\alpha-(d_i-\alpha)\lambda}=1.$$
Technically, it suffices to satisfy only one of these in most cases. The only times where it is necessary to check both is if $\lambda=m_1$ or $\lambda=m_3$. In particular, combining these gives us
\begin{equation}\label{equalszero}
\sum_{i\in[n]}\frac{(d_i-m_1)(d_i-\alpha)}{m_3d_i-m_1\alpha-(d_i-\alpha)\lambda}=0    
\end{equation}
which shows that there can be at most one such $\lambda$ for each $\alpha$ since each term in this sum is monotonic in $\lambda$. We can therefore implicitly define $\lambda$ as a smooth function of $\alpha\in(d_{i_*},d_n)$. 

We see now that for each $\alpha\in(d_{i_*},d_n)$, there exists a unique critical point of our constrained minimization problem, and since $f$ is bounded below by 1 but is not bounded from above, we know this critical point corresponds to a constrained minimum. So we have a smooth curve $\gamma$ parameterized by $\alpha\in(d_{i_*},d_n)$ along which we find the minimum value of $f$ given that 
$\frac{\sum_{i\in[n]}d_ix_i^2}{\sum_{i\in[n]}x_i^2}=\alpha$. We now minimize 
\begin{align*}
f\circ\gamma(\alpha)&=\left(\frac{\alpha}{m_3}\right)^{m_3}\left(\sum_{i\in[n]}\frac{d_i-m_1}{m_3d_i-m_1\alpha-(d_i-\alpha)\lambda}\right)^{m_2}\prod_{i\in[n]}\left(\frac{m_3d_i-m_1\alpha-(d_i-\alpha)\lambda}{d_i-m_1}\right)^{d_i-m_1} \\ 
&=m_3^{-m_3}\alpha^{m_1}\prod_{i\in[n]}\left(\frac{m_3d_i-m_1\alpha-(d_i-\alpha)\lambda}{d_i-m_1}\right)^{d_i-m_1} 
\end{align*}
over $\alpha\in[d_{i_*},d_n]$ to find the unconstrained global infimum. By (\ref{equalszero}), we have 
\begin{align*}
\frac{\mathrm{d}f\circ\gamma}{\mathrm{d}\alpha}&=\frac{m_1}{\alpha}-\sum_{i\in[n]}\frac{(d_i-m_1)\big(m_1-\lambda+(d_i-\alpha)\lambda'\big)}{m_3d_i-m_1\alpha-(d_i-\alpha)\lambda} \\ 
&=\frac{m_1}{\alpha}-\frac{m_1-\lambda}{\alpha}=\frac{\lambda}{\alpha}. 
\end{align*}
So we have a critical point wherever $\lambda(\alpha)=0$. Since 
$\frac{\mathrm{d}^2f\circ\gamma}{\mathrm{d}^2\alpha}=\frac{\lambda'}{\alpha}-\frac{\lambda}{\alpha^2}$
and $\alpha>0$, this critical point is a minimum if and only if $\lambda'>0$ at this point. By differentiating (\ref{equalszero}), we end up with 
$$\sum_{i\in[n]}\frac{(d_i-m_1)\big((d_i-\alpha)^2\lambda'-m_2d_i\big)}{(m_3d_i-m_1\alpha-(d_i-\alpha)\lambda)^2}=0.$$
Solving for $\lambda'$, we have 
$$\lambda'=\frac{\sum_{i\in[n]}\frac{m_2(d_i-m_1)d_i}{(m_3d_i-m_1\alpha-(d_i-\alpha)\lambda)^2}}{\sum_{i\in[n]}\frac{(d_i-m_1)(d_i-\alpha)^2}{(m_3d_i-m_1\alpha-(d_i-\alpha)\lambda)^2}}>0.$$
Thus the global minimum occurs when $\alpha$ satisfies $\lambda(\alpha)=0$, and the $\alpha$ that satisfies this is the unique solution to $\sum_{i\in[n]}\frac{d_i-m_1}{m_3d_i-m_1\alpha}=\alpha^{-1}$. The bound in the statement of the corollary is now achieved by rescaling $\alpha$ by a factor of $m_3$. 
\end{proof}

\newpage 

\printbibliography

\end{document}